\documentclass[a4paper,12pt]{article}
\usepackage[english]{babel}
\usepackage{mathrsfs}
\usepackage{amsfonts}
\usepackage{amsthm}
\usepackage{latexsym,amsmath,amssymb}
\newtheorem{thm}{Theorem}
\newtheorem{prop}[thm]{Proposition}
\newtheorem{lemma}[thm]{Lemma}

\newtheoremstyle{obs}
  {3pt}
  {3pt}
  {}
  {}
  {\bfseries}
  {.}
  {.5em}
  {}
\theoremstyle{obs}
\newtheorem{remark}[thm]{Remark}

\newtheorem{example}[thm]{Example}

\def\qed{\ifvmode\removelastskip\fi
{\unskip\nobreak\hfil\penalty50\hbox{}\nobreak\hfil \hbox{\vrule
height1.2ex width1.2ex}\parfillskip=0pt \finalhyphendemerits=0
\par \smallskip}}

\newcommand{\R}{\mathbb{R}}
\newcommand{\mcR}{\mathcal{R}}
\newcommand{\Ad}{\mathrm{Ad}}

\newcommand{\bu}{{\bar u}}
\newcommand{\bU}{{\bar U}}
\newcommand{\bv}{{\bar v}}
\newcommand{\bx}{\mathbf{x}}
\newcommand{\ue}{u^\epsilon}
\newcommand{\uiae}{u_i^\epsilon}

\newcommand{\spn}{\mathrm{span}}
\newcommand{\codim}{\mathrm{codim\,}}
\newcommand{\Ind}{\mathrm{Ind}}

\title{{\bf Second order conditions for optimality \\ and local controllability of discrete-time systems}}
\author{M. Barbero-Li\~n\'an
\thanks{Departamento de Matem\'aticas, Universidad Carlos III de Madrid,
Avenida de la Universidad 30, 28911 Legan\'es, Madrid, Spain, and Instituto de Ciencias Matem\'aticas (CSIC-UAM-UC3M-UCM). E-mail: mbarbero@icmat.es.}\
 \and
 B. Jakubczyk
 \thanks{Institute of Mathematics, Polish Acadademy of Sciences, \'Sniadeckich 8, 00-956 Warsaw, Poland. E-mail: b.jakubczyk@impan.pl.
\newline{
 Work partially supported by Ministry of Research and Higher Education (Poland), grant N201 607540; by the MICINN (Spain)
MTM2010-21186-C02-02; 2009SGR1338 from the Catalan government; the European project IRSES-project “GeoMech-246981” and the ICMAT Severo Ochoa project SEV-2011-0087.}
\newline{Date: January 25, 2013}}}

\date{}

\begin{document}

\maketitle

\begin{abstract}
We study local controllability and optimal control problems for invertible discrete-time systems.
We present second order necessary and sufficient conditions for optimality and for local controllability.
The conditions are stated in geometric terms, using vector fields naturally associated to the system.
The Hessian of the optimal problem is computed in terms of Lie brackets of vector fields of the system.
\end{abstract}

\vskip-2cm
\section{Introduction}

Nonlinear discrete-time control systems
\[
\Sigma:\ \ \ x(t)=f(x(t-1),u(t))
\]
are much less understood, comparing to continuous-time systems. This is due to the fact that algebraic and geometric tools available in the continuous-time case are not present or, at least, have not been used much for their analysis. Here we have in mind vector fields, Lie bracket and Lie algebraic techniques which are very helpful in the theory of continuous-time systems.

There is a class of discrete-time systems where such tools are available, however. These are \emph{invertible systems},
where the state equations are solvable backwards in time, that is $x(t-1)$ is uniquely defined for given $x(t)$ and $u(t)$.
Families of vector fields can be assigned to such systems (see \cite{89JakSonInBook,90JakSon}) and can be used for analyzing their general controllability properties \cite{93AlbSon,90JakSon}. Invertible systems appear e.g. in modeling of continuous-time systems, as in the case of sampling, and in numerical schemes approximating continuous-time systems. Codistributions can be also assigned to such systems and have been useful for characterizing their observability \cite{2002AlbertiniDAlessandro}.
Recently, differential geometric tools have been used to analyze the accessibility of implicit discrete-time
systems \cite{RiegerSchlacher}.

For continuous-time optimality of a trajectory is closely related to its local controllability
property. Namely, for most optimal control problems, optimal trajectories lie on the boundary of the reachable set.
On the other hand, it is a hard problem to find criteria in terms of the system vector fields which characterize local
controllability along a trajectory (i.e., when the trajectory lies in the interior of the reachable set).
Deep and far reaching results in this direction were obtained in the last forty years
(see e.g.~\cite{Stefani,Hermes,Kawski,Sussmann}) but a complete characterization seems beyond the reach~\cite{99Agrachev,88Sontag}. General second order conditions for optimality in continuous-time can be found in~\cite{Agrachev}.

In this paper we analyze invertible discrete-time systems using the formalism of vector fields from~\cite{90JakSon}.
We present (Section \ref{s-Local controllability}) second order sufficient and necessary conditions for local controllability in terms of those vector fields and their Lie brackets.
In Section~\ref{s-Proof of local controllability} we present basic lemmata which are then used for proving the sufficiency result
in Section~\ref{s-Local controllability}. Second order optimality conditions are presented and proved in Sections~\ref{s-Optimal control}
and~\ref{s-Optimal Hamiltonian}, using the same lemmata. The proof of one version of these results is based on a geometric lemma
in \cite{Agrachev} on local openness of a nonlinear map at a singular point.
We include an illustrating example in Section~\ref{s-Example}.

Our infinitesimal analysis of local controllability, with the use of Lie bracket, can be considered as a starting point to identifying higher order sufficient or necessary conditions for local controllability of discrete-time systems (analogous to conditions in the continuous-time case), in particular to identifying so called ``bad brackets" (see~\cite{Stefani,Hermes,Kawski,97LewisMurray,Sussmann} for the case of continuous-time).

\section{Preliminaries}\label{s-Preliminaries}

Let $M$ and $U$ be two sets called \emph{state space} and \emph{control space}, respectively. A map $f\colon M\times U \rightarrow M$ defines a nonlinear discrete-time control system with the dynamics
\[
 \Sigma:\ \ \ x(t)=f(x(t-1),u(t)),
\]
where $x$ and $u$ are called state and control, respectively, and $t\in \mathbb{Z}$.

We will assume that:
\begin{itemize}
\item[(A1)]
$M$ is an open subset of $\R^n$ or a smooth differentiable manifold of dimension $n$;
\item[(A2)]
$U$ is a subset of $\R^m$, with nonempty interior, and the closure of the interior of $U$ contains $U$;
\item[(A3)]
the map $f\colon M\times U \rightarrow M$ is of class $C^2$;
\item[(A4)]
the system is \emph{invertible}, which means that the maps $f_u:M\to M$, $u\in U$, are diffeomorphisms onto open images.
\end{itemize}

Condition (A3) means that $f$ has a $C^2$ extension to $M\times \widetilde U$, where $\widetilde U\subset \R^m$ is an open superset of $U$. Above, and in the sequel, we denote
\[
 f_u(x)=f(x,u).
\]
The map $f_u:M\to M$ defines the one-step transition defined by control $u$. The invertibility property (A4) means that
$f_u(M)\subset M$ is open and each $f_u:M\to f(M)$ is a $C^2$ diffeomorphism.

Assumption (A4) is needed for associating natural vector fields to system $\Sigma$.
Sometimes we will assume that the system is \emph{strongly invertible}, which means that the maps $f_u:M\to M$, $u\in U$, are diffeomorphisms onto $M$.

It will be convenient to use the notation $t=i$ and write the system equations in the form
\[
 \Sigma:\ \ \ x_i=f_{u_i}(x_{i-1}), \quad \mathrm{where}\ \  x_i=x(i),\ u_i=u(i).
\]
Given an initial state $x_0$ and a control sequence $u_1,\dots,u_N$, the trajectory of $\Sigma$ is defined by the sequence of states $x_1,\dots,x_N$, where $x_i$ is given by the composition of maps $f_{u_i}\circ\cdots\circ f_{u_2}\circ f_{u_1}$ applied to $x_0$. We will usually omit the composition sign and write
$x_i=f_{u_i} \cdots f_{u_1}(x_0)$.
The set of points reachable from $x_0$ in $N$ forward steps is denoted by
\begin{equation*}
 {\mathcal R}^+(x_0,N)=\{x\in M \; | \; \exists \; (u_1,\dots,u_N) \in U^N \; {\rm such \;that} \; x=f_{u_N}\dots f_{u_1}(x_0)\}.
\end{equation*}
Given a subset $\bU\subset U^N$ of control sequences $(u_1,\dots,u_N)$, we define the corresponding reachable set
\begin{equation*}
 {\mathcal R}^+(x_0,\bU)=\{x\in M \; | \; \exists \; (u_1,\dots,u_N) \in \bU \; {\rm such \;that} \; x=f_{u_N}\dots f_{u_1}(x_0)\}.
\end{equation*}

With the purpose to discuss local controllability of discrete-time systems we will
provide an infinitesimal description of local deformations of admissible trajectories of $\Sigma$. We will use vector fields introduced in~\cite{90JakSon} (see also~\cite{81FliessNor,JakNor} for earlier work), where they were needed for characterizing accessibility and controllability.

Assume first that the control is scalar, that it is $U\subset \R$. Then
the following vector fields depending on $u$ can be associated to the invertible system $\Sigma$ (cf.~\cite{90JakSon}),
\begin{eqnarray}
 X^+_u(x)=\left. \dfrac{\partial}{\partial \epsilon} \right|_{\epsilon=0} f_u^{-1} \circ f_{u+\epsilon}(x),&&
 Y^+_u(x)=\left. \dfrac{\partial}{\partial \epsilon} \right|_{\epsilon=0} f_{u+\epsilon}^{-1} \circ f_{u}(x),\label{eq:VFforw} \\
 X^-_u(x)=\left. \dfrac{\partial}{\partial \epsilon} \right|_{\epsilon=0} f_u \circ f_{u+\epsilon}^{-1}(x),\ \ &&
 Y^-_u(x)=\left. \dfrac{\partial}{\partial \epsilon} \right|_{\epsilon=0} f_{u+\epsilon} \circ f_{u}^{-1}(x). \label{eq:VFback}
\end{eqnarray}
The vector fields in~\eqref{eq:VFforw} can be used for the infinitesimal analysis of the variations of forward trajectories of $\Sigma$, whereas the vector fields in~\eqref{eq:VFback} play an analogous role for backward trajectories. Note that $X^+_u$ and $Y^+_u$ are well defined for all $x\in M$, if $\Sigma$ is invertible (for small $\epsilon$ the point $f_{u+\epsilon}(x)$ is in the domain of $f_u^{-1}$, by the implicit function theorem), while the same holds for $X^-_u$ and $Y^-_u$ only when $\Sigma$ is strongly invertible. Due to assumptions (A2) and (A3), these vector fields are well defined for any $u\in U$ and are of class $C^1$.

We will mainly use the $u$-depending vector fields $X^+_u$, however similar results can be obtained to describe local controllability
and optimality by means of the vector fields in~\eqref{eq:VFback} when reachability is defined by
backward trajectories. The vector fields $X^+_u$  can be alternatively defined as
\begin{equation*} \label{eq:X+uAlt}
 X^+_u(x)=\left({\rm d}f_u(x)\right)^{-1} \dfrac{\partial }{\partial u} f_u(x).
\end{equation*}
In the case of multidimensional control these $u$-depending vector fields depend, additionally, on the index $r$ of the component of $u=(u^1,\dots,u^m)\in U$,
\begin{equation*} \label{eq:X+ujAlt}
 X^+_{u,r}(x)=\left({\rm d}f_u(x)\right)^{-1} \dfrac{\partial }{\partial u^r} f_u(x)=\left. \dfrac{\partial}{\partial \epsilon} \right|_{\epsilon=0} f_u^{-1} \circ f_{u+\epsilon e_r}(x),
\end{equation*}
where $e_r$ is the r-th versor in $\R^m$.

Given a vector field $Y$ and a control $u\in U$, we define another vector field using the diffeomorphism
$f_u$ on $M$,
\begin{equation*}\label{Eq:AduY}
 ({\rm Ad}_uY)(x)=({\rm d}f_u(x))^{-1}Y(f_u(x)),
\end{equation*}
where ${\rm d}f_u(x)$ is the differential of $f_u$ evaluated at the point $x$.
(The above definition of Ad, used throughout the paper, is convenient for analyzing forward trajectories. Note that it does not match the usual notation of Lie group theory adopted to left actions.)
More generally, denoting $f_{u_k\dots u_1}(x)=f_{u_k}\dots f_{u_1}(x)$ we define the following vector fields
\begin{equation*}
 ({\rm Ad}_{u_k\dots u_1}Y)(x)=({\rm Ad}_{u_1}\cdots {\rm Ad}_{u_k}Y)(x)= ({\rm d}f_{u_k\dots u_1}(x))^{-1}Y(f_{u_k\dots u_1}(x)).
\end{equation*}
Note that, in the case of scalar control,
\begin{equation*}
  ({\rm Ad}_{u_k\dots u_1}X^+_{u_0})(x)=\left. \dfrac{\partial}{\partial \epsilon} \right|_{ \epsilon=0} f_{u_k\dots u_1}^{-1}
f_{u_0}^{-1}f_{u_0+ \epsilon} f_{u_k\dots u_1} (x).
\end{equation*}
For multidimensional control
\begin{equation*}
  ({\rm Ad}_{u_k\dots u_1}X^+_{u_0,r})(x)=\left. \dfrac{\partial}{\partial \epsilon} \right|_{ \epsilon=0} f_{u_k\dots u_1}^{-1}
f_{u_0}^{-1} f_{u_0+ \epsilon e_r} f_{u_k\dots u_1} (x),
\end{equation*}
where $e_r$ is the r-th versor in $\R^m$.

\begin{prop}{\cite[Proposition 3.2]{90JakSon}} For scalar control the following equalities hold for each $u\in U$:
\begin{enumerate}
\item[(1)]\ \  $X^+_u=-Y^+_u$,\ \ \ \ \ \ \ \  $X^-_u=-Y^-_u$,
\item[(2)]\ \ $X^+_u=-{\rm Ad}_uX^-_u$, \ \ $Y^+_u=-{\rm Ad}_uY^-_u$.
\end{enumerate}
For multidimensional control analogous equalities hold for the vector fields $X^+_{u,r}$, $Y^+_{u,r}$, $X^-_{u,r}$ and $Y^-_{u,r}$.
\end{prop}
For so defined vector fields we can compute their Lie bracket $[\cdot, \cdot]$ in the usual way. We
also denote the Lie bracket of $Y$ and $Z$ as ${\rm ad}Y(Z)=[Y,Z]$.

\begin{prop}{\cite[Proposition 3.3]{90JakSon}}\label{p-Adad} The following equalities hold for any vector field $Z$ and any scalar $u\in U$:
\begin{equation*}
 \dfrac{\partial}{\partial u}{\rm Ad}_uZ={\rm ad}X^+_u({\rm Ad}_uZ),\qquad
 \dfrac{\partial}{\partial u}{\rm Ad}_u^{-1}Z={\rm ad}X^-_u({\rm Ad}_u^{-1}Z). \label{Eq:PartialUAd}
\end{equation*}
In the multidimensional control case these equalities take the form
\begin{equation*}
 \dfrac{\partial}{\partial u^r}{\rm Ad}_uZ={\rm ad}X^+_{u,r}({\rm Ad}_uZ),\qquad
 \dfrac{\partial}{\partial u^r}{\rm Ad}_u^{-1}Z={\rm ad}X^-_{u,r}({\rm Ad}_u^{-1}Z). \label{Eq:PartialUAdmulti}
\end{equation*}
\end{prop}

\section{Local controllability and geometric optimality}\label{s-Local controllability}

We first assume that the control is scalar, $u\in U\subset \R$.
Consider an admissible control sequence $\bu=(u_1,\dots,u_N)$. Given an initial point $x_0$, we say that system $\Sigma$ is \emph{N-step locally controllable} from $x_0$ along the trajectory $x(t)$ corresponding to $\bu$ (shortly, $(x_0,\bu)$-\emph{locally controllable}) if
\begin{equation*}\label{eq-xN1}
x(N)\in \mathrm{int}{\mathcal R}^+(x_0,N).
\end{equation*}
The system
$\Sigma$ will be called  \emph{strongly $(x_0,\bu)$-locally controllable} if for any neighborhood $\bU\subset U^N$ of $\bu$ we have
\begin{equation}\label{eq-xN2}
x(N)\in \mathrm{int}{\mathcal R}^+(x_0,\bU).
\end{equation}
Then $\Sigma$ has this property, with the same $x_0$, for any control sequence $\tilde u$ of length $N'>N$ with $\bu$ being its initial part.

In order to state a sufficient condition for local controllability we need the following notation.
For a fixed control sequence $\bu=(u_1,\dots,u_N)$ we introduce the \emph{first variation vector fields}
\begin{equation*}
Y^i_\bu=Y^i_{u_1\cdots u_i}:=\Ad_{u_1}\cdots\Ad_{u_{i-1}}X^+_{u_i}, \quad i=1,\dots,N
\end{equation*}
(in particular, $Y^1_\bu=X^+_{u_1}$). Note that they have the recurrence property
\begin{equation*}
Y^i_{u_1\cdots u_i}={\rm Ad}_{u_1}Y^{i-1}_{u_2\cdots u_i}.
\end{equation*}
For a given $x$ we introduce the space of vectors in $T_xM$
\begin{equation}
L_\bu(x)=\spn\{Y^i_\bu(x), \quad i=1,\dots,N\}. \label{eq:Lux}
\end{equation}
The family of vector fields defining $L_\bu(x)$ does not necessarily describe a minimal set of generators for such a subspace.
We define a subspace, called the \emph{kernel}, of the vector space of coefficients $a=(a_1,\dots,a_N)$,
\begin{equation*}
K_\bu(x)=\left\{a\in\R^N:\sum_{i=1}^Na_iY^i_\bu(x)=0\right\}.
\end{equation*}
For $i,j\in\{1,\dots,N\}$ we define the \emph{second variation vector fields}
\begin{equation}
Z^{ij}_\bu=Z^{ji}_\bu=\dfrac{1}{2}\,[Y^i_\bu,Y^j_\bu], \qquad i<j, \label{Eq:IIvariation1}
\end{equation}
\begin{equation}
Z^{ii}_\bu=\frac{\partial}{\partial u_i}Y^i_\bu, \ \qquad \label{Eq:IIvariation2}
\end{equation}
where the square bracket denotes the Lie bracket of vector fields on $M$. Equivalently, $Z^{ij}_\bu$ are given by
\[
Z^{ij}_\bu=Z^{ji}_\bu=\dfrac{1}{2}\,\Ad_{u_1}\cdots\Ad_{u_{i-1}}[X^+_{u_i},\Ad_{u_i}\cdots\Ad_{u_{j-1}}X^+_{u_j}], \quad i<j.
\]
Given a point $x_0\in M$, consider the vector-valued quadratic form on the space of parameters
\begin{equation}\label{eq:Hscalar}
H_\bu(a)=H(a)=\sum_{i,j=1}^Na_ia_jZ^{ij}_\bu(x_0)
\end{equation}
(the subscript $\bu$ will be omitted). If $\lambda$ is a covector in $T^*_{x_0}M$, the formula $\lambda\, H(a)=\sum_{i,j}a_ia_j\lambda(Z^{ij}_\bu(x_0))$ defines a real-valued quadratic form.
For a quadratic form $Q$ on a finite dimensional real vector space $V$ we denote by $\Ind^+Q$ (resp. $\Ind^-Q$) the maximal dimension of a subspace $W\subset V$ such that $Q$ restricted to $W$ is strictly positive definite (resp. strictly negative definite). Recall also that $Q$ is indefinite if there are vectors $v$ and $w$ such that $Q(v)>0$ and $Q(w)<0$.

Given a subspace $L\subset T_xM$ we denote by $L^\perp$ its annihilator, $L^\perp=\{\lambda\in T_x^*M\,:\,\lambda\vert_L=0\}$. We will often use a covector $\lambda\in T_{x_0}^*M$ annihilating all vector fields $Y^i_\bu$ at $x_0$, i.e.,
\[
\lambda\, Y^i_\bu(x_0)=0, \quad i=1,\dots,N.
\]
This condition is shortly written as $\lambda\in (L_\bu(x_0))^\perp$.

\begin{thm}\label{t1}
Assume that $\Sigma$ satisfies (A1)-(A4). Given a fixed initial point $x_0\in M$, consider an admissible control sequence
$\bu=(u_1,\dots,u_N)$ such that $\bu\in(\mathrm{int\,} U)^N$ and let $k=\mathrm{codim}\,L_\bu(x_0)$. Then the following statements hold.

\noindent
(a) If $k=1$ and $\lambda\,H$ restricted to $K_\bu(x_0)$  is indefinite, for any $\lambda\in (L_\bu(x_0))^\perp$,  $\lambda\ne0$, then system $\Sigma$ is strongly $(x_0,\bu)$-locally controllable.

\noindent
(b) In general, if $\Sigma$ satisfies the condition
\begin{equation}\label{index condition1}
\Ind^-(\lambda\ H)\vert_{K_\bu(x_0)}\ge k, \quad \forall\ \lambda \in (L_\bu(x_0))^\perp, \ \lambda\ne0,
\end{equation}
then it is strongly $(x_0,\bu)$-locally controllable.
\end{thm}

Note that replacing $\lambda$ with $-\lambda$ gives the same inequality for $\Ind^+$.
Thus condition \eqref{index condition1} says that the quadratic form $\lambda\, H$ on $\R^N$, defined by the second variation vectors at $x_0$, has at least $k$ positive and $k$ negative ``eigenvalues", when restricted to the subspace $K_\bu(x_0)\subset \R^N$.

We will now state a similar result for the multidimensional control, i.e., for $U\subset\R^m$.
We fix a control sequence $\bu=(u_1,\dots,u_N)\in U^N$, where $u_i=(u_i^1,\dots,u_i^m)$. Analogously to the scalar control case we define
\emph{first variation vector fields}
\begin{equation*}
Y^{ir}_\bu=\Ad_{u_1}\cdots\Ad_{u_{i-1}}X^+_{u_i,r},
\end{equation*}
for $i=1,\dots,N$, $r=1,\dots,m$. Given a point $x$ we introduce the space of vectors in $T_xM$
\begin{equation*}
L_\bu(x)=\spn\{Y^{ir}_\bu(x), \quad i=1,\dots,N,\ r=1,\dots,m\}.
\end{equation*}
Consider the space of coefficients $a=(a_1,\dots,a_N)$, where $a_i=(a_i^1,\dots,a_i^m)$. We will use its subspace
\begin{equation*}
K_\bu(x)=\left\{a\in\R^{mN}:\sum_{i=1}^N\sum_{r=1}^m a_i^rY^{ir}_\bu(x)=0\right\},
\end{equation*}
called the \emph{kernel}. For $i,j\in \{1,\dots, N\}$, $r,s\in \{1, \dots, m\}$, we define \emph{second variation vector fields}:
\begin{equation*}
Z^{ir,js}_\bu=Z^{js,ir}_\bu=\dfrac{1}{2}\,[Y^{ir}_\bu,Y^{js}_\bu]
=\dfrac{1}{2}\,\Ad_{u_1}\cdots\Ad_{u_{i-1}}[X^+_{u_i,r},\Ad_{u_i}\cdots\Ad_{u_{j-1}}X^+_{u_j,s}], \quad i<j,
\end{equation*}
\begin{equation*}
Z^{ir,is}_\bu=\Ad_{u_1}\cdots\Ad_{u_{i-1}}\frac{\partial}{\partial u^r_i}X^+_{u_i,s}.
\end{equation*}
For $x_0\in M$ consider the following vector-valued quadratic form on $\R^{Nm}$
\begin{equation}\label{eq:Hmulti}
H_\bu(a)=H(a)=\sum_{ \substack{i,j=1,\dots,N\\ r,s=1,\dots,m}}a^r_ia^s_jZ^{ir,js}_\bu(x_0).
\end{equation}
With the above definitions we have the following result.

\begin{thm}\label{t2}
Theorem \ref{t1} remains valid in the case of multidimensional control.
\end{thm}

The following converse result will be proved in Section  \ref{s-Optimal control} using Theorem \ref{t5}.

\begin{thm}\label{t2'}
If there exists $\lambda\in (L_\bu(x_0))^\perp$ such that $\lambda\,H\vert_{K_\bu(x_0)}$ is positive definite,
\[
\lambda\,H\vert_{K_\bu(x_0)}\ > \ 0,
\]
then $\Sigma$ is not strongly $(x_0,\bu)$-locally controllable.
\end{thm}

Theorems \ref{t1} and \ref{t2} can be used for obtaining necessary conditions on geometric optimality. Recall that, given an initial point $x_0$ and  control sequence $\bu=(u_1,\dots,u_N)$, the corresponding trajectory $x(i)$, $i=0,\dots,N$, of system $\Sigma$ is called \emph{geometrically optimal} if it lies on the boundary of the reachable set, i.e.,
\[
x(i)\in \partial \mcR^+(x_0,i), \quad i=1,\dots,N.
\]
It follows from the invertibility of the system that condition $x(j)\in \mathrm{int}\,\mcR^+(x_0,j)$ implies the inclusion $x(k)\in \mathrm{int}\,\mcR^+(x_0,k)$, for any $k>j$. Thus, the above condition for geometric optimality is equivalent to $x(N)\in \partial \mcR^+(x_0,N)$. Theorems \ref{t1} and \ref{t2} trivially imply

\begin{thm}\label{t3}
If the trajectory of system $\Sigma$ corresponding to initial point $x_0\in M$ and an admissible control sequence
$\bu=(u_1,\dots,u_N)$ such that $\bu\in(\mathrm{int\,} U)^N$ is geometrically optimal, then there exists a nonzero covector $\lambda \in (L_\bu(x_0))^\perp$ such that
\begin{equation*}\label{index condition negative}
\Ind^-(\lambda\ H)\vert_{K_\bu(x_0)}< \mathrm{codim}\,L_\bu(x_0).
\end{equation*}
In particular, if $L_\bu(x_0)$ is the whole tangent space then the trajectory is not geometrically optimal. If $\mathrm{codim}\,L_\bu(x_0)=1$
then it is necessary for geometric optimality that there exists a nonzero $\lambda \in (L_\bu(x_0))^\perp$ (unique up to a positive multiplier) such that the quadratic form $\lambda\ H$ restricted to the kernel $K_\bu(x_0)$ is non-negative definite.
\end{thm}

\newpage

\section{Example}\label{s-Example}

Take $M=\mathbb{R}^3$ and let $U\subseteq \mathbb{R}$ be an open subset. Consider the discrete-time control system on $M$
\begin{equation}\label{eq-system-Ex}
f_u(x,y,z) = \begin{pmatrix} -x+z+\dfrac{u^2}{2}\\ xz-y\\ z+\dfrac{u^2}{2} \end{pmatrix}.
\end{equation}
Denote the state $\bx=(x,y,z)$. As the state is 3-dimensional and the control is scalar, the system is never locally controllable in two steps. It is 3-steps $(\bx,\bu)$-locally controllable if $Y^1_\bu,Y^2_\bu,Y^3_\bu$ are linearly independent at $\bx$ (this simply follows from the inverse function theorem). If they are not linearly independent, the index condition from Theorem \ref{t1} does not help as it never holds here for 3-steps controls. The reader may analyze the cases where Theorem \ref{t2} is applicable.

We will consider 4-steps controls $\bar{u}=(u_1,u_2,u_3,u_4)$. Given an initial condition $\bx=(x,y,z)$,
we have
\begin{equation*}
 {\rm d}f_u(\bx)=\begin{pmatrix} -1 & 0 & 1 \\ z & -1 & x \\ 0 & 0 & 1 \end{pmatrix}, \quad
 {\rm d}f_u^{-1}(\bx)=\begin{pmatrix} -1 & 0 & 1 \\ -z & -1 & x+z \\ 0 & 0 & 1  \end{pmatrix}, \quad
 \dfrac{\partial f}{\partial u}=\begin{pmatrix} u\\ 0 \\ u \end{pmatrix}.
\end{equation*}
The first variation vector fields at $\bx$ are:
\begin{eqnarray*}
 Y^1_\bu(\bx) &=& Y^1_{u_1}(\bx)       = X^+_{u_1}(\bx)= \begin{pmatrix} 0\\ x \\ 1 \end{pmatrix}u_1,\\
 Y^2_\bu(\bx) &=& Y^2_{u_1u_2}(\bx)={\rm Ad}_{u_1} Y^1_{u_2}(\bx)\\
              &=& {\rm d}f_{u_1}^{-1}(\bx)Y^1_{u_2}(f_{u_1}(\bx))=\begin{pmatrix} 1\\2x-\dfrac{1}{2} \,u_1^2\\1 \end{pmatrix}u_2,\\
 Y^3_\bu(\bx) &=& Y^3_{u_1u_2u_3}(\bx)={\rm Ad}_{u_1}\,Y^2_{u_2u_3} (\bx)\\
              &=&{\rm d}f_{u_1}^{-1}(\bx)Y^2_{u_2u_3}(f_{u_1}(\bx))
                                       =\begin{pmatrix} 0 \\ 3x-2z+\dfrac{1}{2}\,u_2^2-u_1^2\\ 1\end{pmatrix}u_3,\\
 Y^4_\bu(\bx) &=& Y^4_{u_1u_2u_3u_4}(\bx)={\rm Ad}_{u_1}\,Y^3_{u_2u_3u_4} (\bx)\\
              &=& {\rm d}f_{u_1}^{-1}(\bx)Y^3_{u_2u_3u_4}(f_{u_1}(\bx))
                                       =\begin{pmatrix} 1 \\ 4x-\dfrac{1}{2}\,u_1^2+u_2^2-\dfrac{1}{2}\, u_3^2\\ 1 \end{pmatrix}u_4.
\end{eqnarray*}
The second variation vector fields can be computed using the definition \eqref{Eq:IIvariation2},
\begin{eqnarray*}
 Z^{11}_\bu(\bx)   &=&  \begin{pmatrix}0 \\ x \\ 1 \end{pmatrix}, \\
 Z^{22}_\bu(\bx)   &=&  \begin{pmatrix} 1 \\  2x -\dfrac{1}{2}u_1^2 \\ 1 \end{pmatrix},\\
 Z^{33}_\bu(\bx)   &=&  \begin{pmatrix}  0 \\ 3x-2z-u_1^2+\dfrac{1}{2}u_2^2\\ 1  \end{pmatrix},\\
 Z^{44}_\bu(\bx)   &=&  \begin{pmatrix} 1 \\ 4x-\dfrac{1}{2}u_1^2+u_2^2-\dfrac{1}{2}u_3^2\\ 1  \end{pmatrix},\\
\end{eqnarray*}
and the definition \eqref{Eq:IIvariation1}: $Z^{ij}_\bu=Z^{ji}_\bu=\frac{1}{2}[Y^i_\bu,Y^j_\bu]$, for $i<j$. Computing the Lie brackets we find that
\[
\begin{array}{ccc}
 Z^{12}_\bu(\bx)   = \dfrac{1}{2}\, \begin{pmatrix}0 \\ -1 \\ 0 \end{pmatrix}u_1u_2, &
 Z^{13}_\bu(\bx)   = \dfrac{1}{2}\, \begin{pmatrix} 0 \\ -2 \\ 0 \end{pmatrix}u_1u_3, &
 Z^{14}_\bu(\bx)   = \dfrac{1}{2}\, \begin{pmatrix} 0 \\ -1 \\ 0 \end{pmatrix}u_1u_4,\\
 Z^{23}_\bu(\bx)   = \dfrac{1}{2}\, \begin{pmatrix} 0 \\ 1 \\ 0 \end{pmatrix}u_2u_3, &
 Z^{24}_\bu(\bx)   = \dfrac{1}{2}\, \begin{pmatrix} 0 \\ 2 \\ 0 \end{pmatrix}u_2u_4,&
 Z^{34}_\bu(\bx)   = \dfrac{1}{2}\, \begin{pmatrix} 0 \\ -1\\ 0 \end{pmatrix}u_3u_4.
\end{array}
\]
{\bf Case I.} Assume that
\begin{equation}\label{eq:assumptionEx1}
u_1\ne0,\ \ u_2\ne0,\ \ u_3\ne0,\ \ u_4\ne0.
\end{equation}
Then the dimension of ${\rm L}_{\bar{u}}(\bx)={\rm span}\left\{ Y^1_{\bu}, Y^2_{\bu}, Y^3_{\bu},Y^4_{\bu}\right\}$ is 3 (then $\Sigma$ is locally controllable) or it is equal to 2, if
\begin{equation} \label{eq:assumptionEx2}
x=\dfrac{1}{4}\, u_3^2-\dfrac{1}{2} \, u_2^2,\quad z=\dfrac{1}{4}\, u_3^2-\dfrac{1}{4}\, u_2^2 - \dfrac{1}{2} u_1^2.
\end{equation}
Let $\bx=(x,y,z)$ satisfy condition \eqref{eq:assumptionEx2}.
Then ${\rm codim}\, {\rm L}_{\bar{u}}(\bx)=1$ and the annihilating space ${\rm L}_{\bar{u}}^\perp(\bx)$ is generated by the covector
\begin{equation}
\lambda=\left(-x+\dfrac{1}{2}u_1^2,1,-x\right)\label{eq:lambdaEx2}.
\end{equation}
The kernel is
\begin{equation*}
 K_{\bu}(\bx)=\left\{ (a_1,a_2,a_3,a_4) \; \colon \; a_3=-\dfrac{u_1}{u_3}\, a_1, \; a_4=-\dfrac{u_2}{u_4} \, a_2 \right\} \subseteq \mathbb{R}^4.
\end{equation*}
For $\lambda$ defined by \eqref{eq:lambdaEx2} the quadratic form $\lambda H$ is given by
\begin{equation*}
\lambda H(a) = \dfrac{1}{2}\,
\begin{pmatrix} a_1 \\ a_2 \\ a_3 \\ a_4\end{pmatrix}^T
\begin{pmatrix}
 0 & -u_1u_2 & -2u_1u_3 & -u_1u_4 \\
 -u_1u_2 & 0 & u_2u_3  & 2u_2u_4 \\
 -2u_1u_3 & u_2 u_3 & A & -u_3 u_4 \\
 -u_1 u_4 & 2u_2u_4 & -u_3u_4 & B \end{pmatrix}\,
\begin{pmatrix} a_1 \\ a_2 \\ a_3 \\ a_4\end{pmatrix},
\end{equation*}
where
\[
\dfrac{1}{2}\,A=2x-2z-u_1^2+\dfrac{1}{2}u_2^2, \qquad \dfrac{1}{2}\,B=2x+u_1^2-\dfrac{1}{2}u_3^2.
\]
Assuming \eqref{eq:assumptionEx1}, we have $A=0=B$ at points $\bx\in M$ satisfying condition \eqref{eq:assumptionEx2}, where $\dim\, L_\bu(\bx)=2$.  Then
\begin{equation*}
\lambda H(a)= \dfrac{1}{2}\,
\begin{pmatrix} A_1 \\ A_2 \\ A_3 \\ A_4\end{pmatrix}^T
\begin{pmatrix}
 0 & -1 & -2 & -1 \\
 -1 & 0 & 1  & 2 \\
 -2 & 1 & 0 & -1 \\
 -1 & 2 & -1 & 0 \end{pmatrix}\,
\begin{pmatrix} A_1 \\ A_2 \\ A_3 \\ A_4\end{pmatrix}, \qquad \mathrm{where} \quad A_i=a_iu_i.
\end{equation*}
From the form of the kernel $K_\bu$ we see that, for fixed $\bu$ satisfying \eqref{eq:assumptionEx1}, the vector $(a_1,a_2,a_3,a_4)$ is in the kernel if and only if $A_1+A_3=0$ and $A_2+A_4=0$. Therefore, taking $A_3=-A_1$ and $A_4=-A_2$ we get
\[
(\lambda H\vert_{K_\bu})(a)= 4(A_1^2-A_1A_2-A_2^2).
\]
This expression is indefinite, treated as a quadratic function of the vector $(a_1,a_2)$ parameterizing the kernel $K_\bu$, which means that there exist values of $a_1$ and $a_2$ such that the quadratic form is positive and negative. Thus, by  Theorem~\ref{t1}, for any 4-step control $\bu$ satisfying \eqref{eq:assumptionEx1} and any initial point fulfilling \eqref{eq:assumptionEx2} the system is strongly $(\bx,\bu)$-locally controllable, even if \eqref{eq:assumptionEx2} means that it does not satisfy the first order sufficient condition for local controllability.

\noindent
{\bf Case IIa.}
For other 4-steps controls, not satisfying \eqref{eq:assumptionEx1}, we consider only one case where $u_1=u_3=0$ and $u_2\neq0$, $u_4\neq0$. In this case $\dim\,L_\bu(\bx)=1$, if the initial condition satisfies $x=-\dfrac{1}{2}u_2^2$, and $\dim\,L_\bu(\bx)=2$, otherwise. If $x\neq-\dfrac{1}{2}u_2^2$, we have $k=\codim\,L_\bu(\bx)=1$. Then $(L_\bu(\bx))^\perp$ is generated
by the covector $\lambda=(-1,0,1)$, the kernel $K_\bu(\bx)$ consists of vectors $(a_1,0,a_3,0)^T$, $a_1,a_3\in\R$, and the quadratic form $\lambda\,H$ evaluated on vectors in $K_\bu(\bx)$ is
\[
\lambda H(a)=a_1^2+a_3^2.
\]
We see that $\lambda H$ restricted to $K_\bu(\bx)$ is positive definite.  It then follows from Theorem \ref{t2'} that the system  is not strongly $(\bx,\bu)$-locally controllable in this case.

\noindent
{\bf Case IIb.} If $u_1=u_3=0$, $u_2\neq0$, $u_4\neq0$ and $x=-\dfrac{1}{2}u_2^2$, Theorem \ref{t2'} can be used again. Namely, the annihilator $L_\bu^\perp(\bx)$ is spanned by the covectors $\lambda_1=(-1,0,1)$, $\lambda_2=(2x,-2,2x)$ and $\lambda\in L_\bu^\perp(\bx)$
has the general form
\[
\lambda=a\lambda_1+b\lambda_2=(-a+2bx,-2b,a+2bx).
\]
The corresponding quadratic form is
\[
\lambda\,H=
\begin{pmatrix}
 a & 0 &  0  & 0 \\
 0 & 0 &  0  & -2b \\
 0 & 0 &a+b(4z+u_2^2)& 0 \\
 0 & -2b & 0 & 0
\end{pmatrix}.
\]
The kernel $K_\bu(\bx)$ is spanned by the vectors
\[
V_1= (1, 0, 0, 0)^T, \quad
V_2= (0, 0, 1, 0)^T, \quad
V_3= (0, 1, 0, -1)^T.
\]
In the basis $V_1,V_2,V_3$ the matrix $\{\lambda\,H(V_i,V_j)\}$ of the quadratic form $\lambda H\vert_{K_\bu(\bx)}$ is
\[
\lambda H\vert_{K_\bu(\bx)}=\mathrm{diag}\{a,a+b(4z+u_2^2),4b\}.
\]
It is positive definite if $a\gg b>0$. Using Theorem \ref{t2'} we see that also in this case the system  is not strongly $(\bx,\bu)$-locally controllable.

{\bf Conclusions.} (a) For any 4-steps control $\bu$ satisfying \eqref{eq:assumptionEx1} (thus, for any N-steps control with initial part $\bu$) and any initial condition $\bx$ the system \eqref{eq-system-Ex} is strongly $(\bx,\bu)$-locally controllable.
\newline
(b) If $\bu$ satisfies $u_1=u_3=0$ and $u_2\ne0\ne u_4$ then \eqref{eq-system-Ex} is not strongly $(\bx,\bu)$-locally controllable.

\section{Proof of local controllability}\label{s-Proof of local controllability}

We will prove the sufficiency result in Theorem~\ref{t1}, only (the proof of Theorem \ref{t2} is analogous).
The proof, as well as further optimality results, are based on the following lemmata.

Given a control sequence $\bu=(u_1,\dots,u_N)$ we define the composed map $f_\bu:M\to M$,
\[
f_\bu=f_{u_N}\cdots f_{u_1},
\]
which is a diffeomorphism. Let $W=\mathrm{int}\,U\times\cdots\times\mathrm{int}\,U\subset \R^{N}$ and consider the map $F:W\to \R^n$ defined for a fixed $x_0\in M$ by
\begin{equation}\label{eq-def-of-F}
F(\bu)=f_\bu(x_0),
\end{equation}

\begin{lemma}\label{le-image} We have
\[ \mathrm{Im}\,dF(\bu)=df_\bu(x_0)L_\bu(x_0). \]
\end{lemma}

\begin{lemma}\label{le-kernel} The kernel of $dF(\bu)$ is given by
\[ \ker dF(\bu)=K_\bu(x_0). \]
\end{lemma}

\begin{lemma}\label{le-Hessian} The second differential of $F$ at $\bu$, restricted to the kernel $\ker dF(\bu)=K_\bu(x_0)$, coincides with $df_\bu(x_0)H$ restricted to this kernel,
\[
d^2F(\bu)\vert_{K_\bu(x_0)} = df_\bu(x_0)H\vert_{K_\bu(x_0)}.
\]
\end{lemma}

{\bf Proof of Lemma \ref{le-image}.}
For $x_0\in M$ and a control sequence $\bu=(u_1,\dots,u_N)\in (\mathrm{int}\,U)^N$  denote $x_i=f_{u_i}\cdots f_{u_1}(x_0)$, $i=1,\dots,N$. Then $x_i=f_{u_i}(x_{i-1})$ and the image of $dF(\bu)$ is spanned by the vectors
\begin{eqnarray*}
\frac{\partial}{\partial u_i}f_{u_N}\cdots f_{u_i}\cdots f_{u_1}(x_0)
 & = & d(f_{u_N}\cdots f_{u_{i+1}})(x_i)\frac{\partial}{\partial u_i} f_{u_i}(f_{u_{i-1}}\cdots f_{u_1}(x_0))\\
 & = & df_\bu(x_0)(df_{u_1}(x_0))^{-1}\cdots(df_{u_i}(x_{i-1}))^{-1}\frac{\partial}{\partial u_i}f_{u_i}(f_{u_{i-1}}\cdots f_{u_1}(x_0))\\
& = & df_\bu(x_0)(\Ad_{u_1}\cdots\Ad_{u_{i-1}}X^+_{u_i})(x_0)=df_\bu(x_0)Y^i_\bu(x_0).
\end{eqnarray*}
This proves the lemma because of the definition of $L_\bu(x_0)$ in~\eqref{eq:Lux}. \qed

\medskip

{\bf Proof of Lemma \ref{le-kernel}.}
Given a control sequence $\bu=(u_1,\dots,u_N)\in (\mathrm{int}\,U)^N$ and a vector of parameters $a=(a_1,\dots,a_N)\in \R^N$, denote
\[
\uiae=u_i+a_i\epsilon, \quad i=1,\dots,N,
\]
where $\epsilon$ is a small real parameter. Fix $x_0\in M$. In order to find the kernel of the differential $dF(\bu)$ it is enough to find such $a=(a_1,\dots,a_N)$ that $\partial f_{\ue_N}\cdots f_{\ue_1}(x_0)/\partial \epsilon=0$, at $\epsilon=0$. We compute
\[
\frac{\partial}{\partial \epsilon}f_{\ue_N}\cdots f_{\ue_1}(x_0)\vert_{\epsilon=0}  =
 \sum_{i=1}^N a_i\, \frac{\partial}{\partial u_i}f_{u_N}\cdots f_{u_i}\cdots f_{u_1}(x_0) = df_\bu(x_0)\left(\sum_{i=1}^N a_i Y^i_\bu(x_0)\right)
\]
where we use the equality established in the preceding proof. Since $df_\bu(x_0)$ is an isomorphism, the above sum is equal to zero if and only if $\sum_{i=1}^N a_i Y^i_\bu(x_0)=0$, which ends the proof. \qed

\medskip

{\bf Proof of Lemma \ref{le-Hessian}.} We will use the notation and the results from the preceding proofs.
In order to compute $d^2F(\bu)$ on the kernel $K_\bu(x_0)$ it is enough to compute the second order derivative  $\frac{\partial^2}{\partial \epsilon^2} f_{\ue_N}\cdots f_{\ue_1}(x_0)\vert_{\epsilon=0}$, for any vector $a=(a_1,\dots,a_N)$ satisfying $\sum_{i=1}^N a_i Y^i_\bu(x_0)=0$. We have
\begin{eqnarray*}
\frac{\partial^2}{\partial \epsilon^2}f_{\ue_N}\cdots f_{\ue_1}(x_0)\vert_{\epsilon=0}
 & = & \sum_{i,j=1}^N a_ia_j\, \frac{\partial}{\partial u_i}\frac{\partial}{\partial u_j}f_{u_N}\cdots f_{u_1}(x_0)\\
 & = &
 \sum_{i=1}^Na_i\frac{\partial}{\partial u_i}\left(\sum_{j=1}^N a_j\, \frac{\partial}{\partial u_j}f_{u_N}\cdots f_{u_1}(x_0)\right)\\
 & = &   \sum_{i=1}^Na_i\frac{\partial}{\partial u_i}\left( df_\bu(x_0)\left(\sum_{j=1}^N a_j Y^j_\bu(x_0)\right)\right)\\
 & = &    df_\bu(x_0) \left( \sum_{i=1}^N \sum_{j=1}^N a_ia_j \frac{\partial}{\partial u_i} Y^j_\bu(x_0)\right),
\end{eqnarray*}
where in the last two equalities we use the equality shown in the proof of Lemma \ref{le-kernel} and the fact that $\sum_j a_j Y^j_\bu(x_0)=0$. If $i<j$ then, using Proposition \ref{p-Adad}, we get
\begin{eqnarray*}
\frac{\partial}{\partial u_i} Y^j_\bu(x_0)
  & = &\frac{\partial}{\partial u_i}\left(\Ad_{u_1}\cdots\Ad_{u_{j-1}}X^+_{u_j}\right)(x_0)\\
  & = & \left(\Ad_{u_1}\cdots\Ad_{u_{i-1}}[X^+_{u_i}, \Ad_{u_i}\cdots\Ad_{u_{j-1}}X^+_{u_j}]\right)(x_0) = [Y^i_\bu,Y^j_\bu](x_0) = 2Z^{ij}_\bu(x_0).
\end{eqnarray*}
For $i=j$ we have
\[
\frac{\partial}{\partial u_i} Y^i_\bu(x_0) = Z^{ii}_\bu(x_0).
\]
Finally, $\frac{\partial}{\partial u_i} Y^j_\bu(x_0)=0$ if $i>j$.
Thus, for $a\in K_\bu(x_0)$, we obtain
\[
(d^2F(\bu))(a)=
\frac{\partial^2}{\partial \epsilon^2}f_{\ue_N}\cdots f_{\ue_1}(x_0)\vert_{\epsilon=0} =  df_\bu(x_0) \left( \sum_{i=1}^N \sum_{j=1}^N a_ia_j Z^{ij}_\bu(x_0)\right)=df_\bu(x_0) H (a)
\]
and the proof is complete. \qed

The proof of Theorem \ref{t1} is based on Lemmata \ref{le-image}, \ref{le-kernel}, \ref{le-Hessian} and the following result \cite[Theorem 20.3]{Agrachev}.

\begin{lemma}\label{le-main}
Let $F:W\to \R^n$ be a map of class $C^2$, where $W\subset\R^k$ is an open subset. Fix a point $w_0\in W$ and denote the kernel and the image of the differential $dF(w_0)$ of $F$ at $w_0$ by
\[
K=\ker\, dF(w_0) \quad \mathrm{and}\ \  L=\mathrm{Im}\,dF(w_0).
\]
If $k:=\codim L\ge 1$ and
\[
\Ind^- \lambda\, d^2F(w_0)\vert_K\ge k \quad \forall\ \lambda \in L^\perp\subset (\R^n)^*,\ \ \lambda\ne0,
\]
then $F(w_0)\in \mathrm{int}\, F(W)$. (Here $\lambda\, d^2F(w_0)\vert_K$ denotes the quadratic form $\lambda\, d^2F(w_0)$ restricted to $K$.)
\end{lemma}

{\bf Proof of Theorem \ref{t1}.}
Consider the map \eqref{eq-def-of-F}. By the definition of strong local controllability~\eqref{eq-xN2} it is enough to show that, for a control sequence $\bu=(u_1,\dots,u_N)$, the point $F(\bu)$ is in the interior of the image $F(W)$ of some neighborhood $W$ of $\bu$.
This is guaranteed if the assumptions of Lemma  \ref{le-main} are satisfied, for $w_0=\bu$.
From the hypothesis of Theorem \ref{t1} it follows that we can use Lemmata \ref{le-image}, \ref{le-kernel}, \ref{le-Hessian}
and we see that the assumption of Lemma \ref{le-main} is fulfilled. This completes the proof.  \qed

\section{Optimal control}\label{s-Optimal control}

The general local characteristics of the system, defined at a given initial state $x_0$ and corresponding to a given control sequence $\bu$ (see Section \ref{s-Local controllability}), may also be used for analyzing optimality of the control. In this case standard tools from optimization theory are available. Choosing simple optimal control problems we indicate the use in these problems of the geometric objects introduced earlier: the first variation vector fields $Y^i_\bu$ and the space $L_\bu(x_0)$ spanned by them at $x_0$, the kernel $K_\bu(x_0)$, the second variation vector fields $Z^{ij}_\bu$ and the corresponding Hessian matrix $H$.

Consider the following optimal control problems. Given a system
\begin{equation*}
 \Sigma:\ \  x(t)=f(x(t-1),u(t)), \quad x(0)=x_0, \ \ x(t)\in M,\ u(t)\in U\subseteq \mathbb{R}^m,
\end{equation*}
satisfying conditions (A1)-(A4) and given an initial point $x_0\in M$, find a control sequence $u(t)$, $t=1,\dots,N$,
which minimizes a function
\[
(P1):\ \ \ \ \ \ \varphi(x(N)). \qquad\qquad\qquad\qquad
\]
The function $\varphi:M\to\R$, called \emph{final cost}, is assumed of class $C^2$. The number of steps $N$ is assumed fixed.
Another version of the problem is obtained if, instead, we minimize a \emph{cost functional}
\[
(P2):\ \ \ \varphi(x(N))+\sum_{t=1}^Nc(x(t-1),u(t)),
\]
where $c:M\times U\to \R$ is called \emph{cost function} and assumed of class $C^2$. With analogy to continuous-time systems we will call (P1) Meyer problem and (P2) Bolza problem.

The Bolza problem can be reduced to the Meyer problem by introducing an additional state coordinate $x^0(t)=\sum_{i=1}^tc(x(i-1),u(i))$. This coordinate satisfies the additional state equation
\begin{equation*}
x^0(t)=x^0(t-1)+c(x(t-1),u(t)), \quad x^0(0)=0, \label{eq:x0}
\end{equation*}
and then problem (P2) is equivalent to problem (P1) with augmented state $\hat x=(x^0,x)$  and the final cost  $\hat \varphi(\hat x(N))=\varphi(x(N))+x^0(N)$.

We will state our results for problem (P1), only, using the notation introduced in Section \ref{s-Local controllability}. As earlier, we denote $u(i)=u_i$, $x(i)=x_i$, and $f_u(x)=f(x,u)$. A control sequence $\bu=(u_1,\dots,u_N)$ is called \emph{locally optimal} for (P1) or (P2) if it is optimal among all sequences in a neighborhood $\bar U\subset U^N$ of $\bu$.

\begin{thm}\label{t4}
If $\bu$ is a locally optimal control sequence for problem (P1) and $\bu\in (\mathrm{int}\,U)^N$ then the covector
\begin{equation}\label{eq-lambda}
\lambda=d\varphi(x_N)\,df_{u_N}(x_{N-1})\cdots\cdots df_{u_1}(x_0)
\end{equation}
satisfies
\[
\lambda\, Y^i_\bu(x_0)=0, \quad i=1,\dots,N, \eqno{(I)}
\]
and
\[
\lambda\, H\vert_{K_\bu(x_0)}\ \ge \ 0. \eqno{(II)}
\]
\end{thm}

In condition (II) the inequality means that the quadratic form is non-negative definite.
Conditions (I) and (II) can be called, respectively, first order and second order necessary conditions for local optimality.
Clearly, (I) is equivalent to $\lambda\in (L_\bu(x_0))^\perp$. A converse result is more complicated.

\begin{thm}\label{t5}
Given an initial point $x_0$ and a control  $\bu\in (\mathrm{int}\,U)^N$, let $\lambda$ be the covector defined by \eqref{eq-lambda}. Assume that (I) holds and condition (II) is strengthened to
\[
 \lambda\, H\vert_{K_\bu(x_0)}\ > \ 0. \eqno{(III)}
\]
Then there exists a quadratic form  $Q$ on the image $L\subset T_{x_N}M$ of the composed map
\[
df_{u_N}(x_{N-1})\cdots\cdots df_{u_1}(x_0)
\]
such that if
\[
d^2\varphi(x_N)\vert_L\ > \ Q \eqno{(IV)}
\]
then $\bu$ is locally optimal for problem (P1). Here $Q$ depends on system $\Sigma$ and $x_0$, $\bu$ and $d\varphi(x_N)$.
\end{thm}

\begin{remark}\label{rem13}
Condition (IV) means that the quadratic form $d^2\varphi(x_N)\vert_L-Q$ is positive definite. The form $Q$ will be determined by formulae \eqref{eq-Q}, \eqref{eq-Q'} in the proof. Note that if (IV) does not hold for a given $\varphi$ then it is satisfied for another $\varphi$ with the same $d\varphi(x_N)$ and suitable $d^2\varphi(x_N)$.
\end{remark}

\begin{remark}
Clearly, the problem (P1) can be reduced to a standard optimization problem. Then standard second order conditions for optimality
can be used, as will be seen at the beginning of the proof of Theorem \ref{t4}. However, these conditions are impractical in use
as they involve multiple compositions of nonlinear maps. Neither they give much geometric insight into the problem.
\end{remark}

{\bf Proof of Theorem \ref{t2'}.}
The assumptions of the theorem imply that there is a covector $\lambda$ at $x_0\in M$ such that conditions (I) and $\lambda H\vert_{K_\bu(x_0)}>0$ hold. Define a covector $\tilde\lambda$ at $x_N=f_{u_N}\cdots f_{u_1}(x_0)$,
\[
\tilde\lambda=\lambda\,(df_{u_1}(x_0))^{-1}\cdots(df_{u_N}(x_{N-1}))^{-1},
\]
and the linear function $\varphi(x)=\tilde\lambda x$. Then  $x_0$, $\bu$ and $\varphi$ satisfy assumptions (I) and (III) in Theorem \ref{t5}. Let $\tilde Q$ be a symmetric matrix such that $\tilde Q\vert_L$ satisfies condition (IV) in Theorem \ref{t5}. Define another function
\[
\tilde \varphi(x)=\tilde\lambda x + \frac{1}{2}(x-x_N)^T(\tilde Q+\varepsilon I)(x-x_N),
\]
where  $I$ is identity matrix and $\varepsilon>0$. Then $\tilde \varphi$ satisfies condition (IV), too, (with $Q$ replaced by $\tilde Q$) and Theorem \ref{t5} implies that control $\bu$ is locally optimal, for problem (P1) with final cost $\tilde\varphi$.

The function $\tilde \varphi$ has regular level sets in a neighborhood of $x_N$. Since optimality of $\bu$ implies that the whole local reachable set from $x_0$ (obtained using controls $\bv$ in a neighborhood of $\bu$) lies above or on the level set of the minimal value, no neighborhood of $x_N$ is covered by the local N-step reachable set. Thus $\Sigma$ is not strongly $(x_0,\bu)$-locally controllable.
\qed

{\bf Proof of Theorem \ref{t4}.}
We use the notation introduced in Section \ref{s-Proof of local controllability}, formula~\eqref{eq-def-of-F}, and assume that a local coordinate system is chosen in a neighborhood of $x_N\in M$. We have
\[
x_N=F(\bu)=f_\bu(x_0),
\]
where $f_\bu$ is the composition $f_\bu=f_{u_N}\cdots f_{u_1}$.
Clearly, problem (P1) is equivalent to minimization of the composed function
\[
\psi(\bu)=\varphi\circ F(\bu).
\]
Since $d\psi(\bu)=d\varphi(x_N)\,dF(\bu)$, the first order necessary condition for minimum can be written as
\begin{equation}\label{eq-Iorder}
d\varphi(x_N)\,dF(\bu)=0
\end{equation}
and the second order condition
\begin{equation}\label{eq-IIorder-full}
d^2\psi(\bu)=d\varphi(x_N)\,d^2F(\bu)+d^2\varphi(x_N)(dF(\bu)\ \cdot\ , dF(\bu)\ \cdot\ )\ \ge\ 0.
\end{equation}
Above we treat $d^2F(\bu)$ as a vector-valued symmetric bilinear form and $d^2\varphi(x_N)$ as a symmetric bilinear form. The inequality means that the corresponding quadratic form is non-negative definite. If the above quadratic form is restricted to the kernel of $dF(\bu)$, equal to the space $K_\bu(x_0)$ by Lemma~\ref{le-kernel}, then the second term vanishes and we get the condition
\begin{equation}\label{eq-IIorder}
d\varphi(x_N)\,d^2F(\bu)\vert_{K_\bu(x_0)}\ \ge\ 0.
\end{equation}
We will show that conditions \eqref{eq-Iorder} and \eqref{eq-IIorder} are equivalent to assertions (I) and (II) of the theorem.

We first prove equivalence of conditions \eqref{eq-Iorder} and (I). Note that \eqref{eq-Iorder} can be written in the form
\begin{equation}\label{eq-I'}
d\varphi(x_N)\, \mathrm{Im}\,dF(\bu)=0.
\end{equation}
Since the covector $\lambda$ in \eqref{eq-lambda} is equal to $\lambda = d\varphi(x_N)\, df_\bu(x_0)$
and $\mathrm{Im}\,dF(\bu)=df_\bu(x_0)L_\bu(x_0)$, by Lemma~\ref{le-image}, we can write \eqref{eq-I'}  as
\[
\lambda\, L_\bu(x_0)=0,
\]
which is equivalent to condition (I).

The equality  $d^2F(\bu)\vert_{K_\bu(x_0)} = df_\bu(x_0)H\vert_{K_\bu(x_0)}$ from Lemma \ref{le-Hessian} implies that condition \eqref{eq-IIorder}
is equivalent to
\begin{equation}\label{eq-IIorderBis}
\lambda\, H\vert_{K_\bu(x_0)}\ \ge \ 0
\end{equation}
which is condition (II) in the theorem.
\qed

{\bf Proof of Theorem \ref{t5}.}
We use the notation from the above proof. Notice that condition \eqref{eq-Iorder} and condition \eqref{eq-IIorder-full} strengthened to strong inequality (positive definiteness) are standard sufficient conditions for local optimality of the point $\bu$, for the function $\psi$. Thus it is enough to show that \eqref{eq-Iorder} and the strengthened version of \eqref{eq-IIorder-full} follow from the assumptions of the theorem.

We have seen in the first part of the proof of Theorem \ref{t4} that assumption (I) is equivalent to condition \eqref{eq-Iorder}.
Thus we should prove that, with appropriately chosen $d^2\varphi(x_N)$, the strengthened version of \eqref{eq-IIorder-full} holds. This means that, for any nonzero vector $v$, we should have
\begin{equation}\label{eq-positiveA+B}
d\varphi(x_N)\,d^2F(\bu)(v,v)+d^2\varphi(x_N)(dF(\bu)v,dF(\bu)v)\ > \ 0.
\end{equation}

In the standard basis in $\R^{mN}$ we can identify the bilinear forms appearing in \eqref{eq-positiveA+B} with symmetric matrices
\[
A\simeq d\varphi(x_N)\,d^2F(\bu), \quad B\simeq d^2\varphi(x_N)(dF(\bu)\,\cdot\, ,dF(\bu)\,\cdot\,).
\]
Denote for brevity $K=K_\bu(x_0)$.
Define the subspace $\ker A=\{v\in \R^{mN}:Av=0\}\subset \R^{mN}$ and notice that $\ker A\cap K=\{0\}$.
This follows from condition (III). Namely, the quadratic form $A$ restricted to $K$ is equal to $\lambda H\vert_K$, by Lemma \ref{le-Hessian} and the definition \eqref{eq-lambda} of $\lambda$. Thus $A\vert_K$ is positive definite. Since $A$ restricted to $\ker A$ is zero, the intersection of $\ker A$ and $K$ must be trivial.

We will bring the bilinear form $A$ to a block-diagonal form.
Since $\ker A\cap K=\{0\}$, we can choose a complement $E$ of $\ker A$ in $\R^{mN}$ so that $K\subset E$.
Since $A$ is nondegenerate on $E$ and it is positive definite  on $K$, the $A$-orthogonal complement $K^\perp=\{v\in E:v^TAw=0\ \forall\ w\in K\}$ of $K$ in $E$ has trivial intersection with $K$. Thus, $E=K\oplus K^\perp$ and $\R^{mN}$ is
the direct sum
\[
\R^{mN}=K\oplus K^\perp\oplus \ker A.
\]
From the definitions of $\ker A$ and $K^\perp$ it follows that, in this  decomposition,
\[
A=\begin{pmatrix} A_{11} & 0 & 0 \\ 0 & A_{22} & 0 \\ 0 & 0 & 0 \end{pmatrix},  \qquad
B=\begin{pmatrix} 0 & 0 & 0 \\ 0 & B_{22} & B_{23} \\ 0 & B_{32} & B_{33} \end{pmatrix},
\]
where the zeros in $B$ appear because $v^TBw=0$ if one of the vectors $v,w$ belongs to $K=\ker dF(\bu)$.

From the definition of $B$ and the fact that $dF(\bu)$ is injective on the complement $K^\perp\oplus \ker A$ of $K=\ker dF(\bu)$
we deduce that the blocks $B_{ij}$ can be chosen arbitrarily, with an appropriate choice of $d^2\varphi(x_N)$.
Our aim is to characterize those $d^2\varphi(x_N)$ which make $B$ satisfying \eqref{eq-positiveA+B}, that is, $A+B>0$.

The weaker condition $A+B\ge 0$ is satisfied, if we choose $B=\tilde B$ with $\tilde B_{22}=-A_{22}$ and the remaining $\tilde B_{ij}=0$, since $A_{11}=A\vert_K=\lambda\,H\vert_K>0$, as proved earlier. Such matrix $\tilde B$ is obtained with choosing $\varphi=\tilde\varphi$ so that $d^2\tilde \varphi(x_N)=\tilde Q$, where $\tilde Q$ is a symmetric matrix characterized via the equality
\begin{equation}\label{eq-Q}
S^T\tilde QS=\tilde B
\end{equation}
and $S$ is the matrix of the linear map $dF(\bu):\R^{mN}\to T_{x_N}M$, in a linear basis in $T_{x_N}M$.

Now we define
\begin{equation}\label{eq-Q'}
Q=\tilde Q\vert_L,
\end{equation}
that is, $Q$ is the quadratic form determined by the matrix $\tilde Q$ in \eqref{eq-Q}, restricted to the subspace $L=\mathrm{Im}\,dF(\bu)=\mathrm{Im}\,S$, and $\varphi$ is chosen so that $d^2\varphi(x_N)\vert_L>Q$. Then we have $A+B>0$, where $B=S^Td^2\varphi(x_N)S$. Indeed, then $S^Td^2\varphi(x_N)S > S^T\tilde QS$ (since the map defined by $S$ is injective) and $A+B$ is block-diagonal, with two positive definite blocks, thus $A+B>0$.
\qed

\begin{remark}
It follows from the above proof that the quadratic form $Q$ defined by \eqref{eq-Q}, \eqref{eq-Q'} has the property that for any other
quadratic form $Q'$ having the property in the theorem we have $Q'\ge Q$. The form $Q$ is uniquely defined by this minimality property.
\end{remark}

\begin{example}
Consider the system given by the dynamics \eqref{eq-system-Ex} and a final cost function $\varphi:M\to\R$.
Then, for any initial condition, a control sequence $\bu=(u_1,u_2,u_3,u_4)$ such that $u_i\ne0$ for all $i$ can not be optimal for problem (P1) if only $d\varphi(x_N)\ne0$ at the final point $x_N$. To prove suppose that $\bu$ is optimal.
Then the covector $\lambda$ given by \eqref{eq-lambda} satisfies condition (I), thus $\lambda\in (L_\bu(\bx))^\perp$.
We have checked in Section \ref{s-Example} (Case I) that either $\dim L_\bu(\bx)=3$ (then $\lambda=0$ and $d\varphi(x_N)=0$),
or $\dim L_\bu(\bx)=2$ and then any $\lambda\in(L_\bu(\bx))^\perp$ is unique up to a multiplier, thus it can be taken as in formula \eqref{eq-lambda}. It was checked in Section \ref{s-Example} that, for such nonzero $\lambda$, the quadratic form $\lambda\,H\vert_{K_\bu(\bx)}$ is indefinite. Thus the necessary condition (III) in Theorem \ref{t4} is not satisfied and we conclude that
$\bu$ is not an optimal control, if $d\varphi(x_N)\ne0$.

Cases IIa and IIb in the example in Section \ref{s-Example} admit optimal control sequences. This can be seen from Theorem \ref{t5}. Namely, in Case IIa we have $\codim\,L_\bu(\bx)=1$ and covectors $\lambda\in(L_\bu(\bx))^\perp$ (equivalently $\lambda$ satisfying condition (I)) are unique, up to multiplier. Let $d\varphi(x_N)$ be such that the corresponding $\lambda$ given by \eqref{eq-lambda} satisfies the first order condition (I). Then  $\lambda\in (L_\bu(\bx))^\perp$ and we have checked that, for such $\lambda$, the quadratic form $\lambda H$ restricted to the kernel is positive definite, thus it satisfies the sufficiency condition (III). From Theorem \ref{t5} we deduce that
the control $\bu$ is optimal, provided that $d^2\varphi(x_N)>\tilde Q$, for some symmetric matrix $\tilde Q$ depending on $\bu$ and the initial state (this matrix is determined from condition \eqref{eq-Q}). In fact, it is enough that $d^2\varphi(x_N)\vert_L>\tilde Q\vert_L$, where $L=\mathrm{Im}dF(\bu)$. The conclusion in Case IIb is similar.
\end{example}

\section{Optimality conditions using Hamiltonian}\label{s-Optimal Hamiltonian}

We will now state second order necessary conditions for optimal control problems applying a formalism similar to Hamiltonian formalism used for continuous-time systems. In order to be able to use duality of vectors and covectors we assume that $M$ is an open subset of an affine or a vector space. In this case the difference $x(t)-x(t-1)$ of two consecutive states can be treated as a vector.

For simplicity we will assume that $M$ is an open subset of a vector space,
\[
M\subset \R^n.
\]
In this case, we can canonically identify all the spaces of (tangent) vectors so that $T_xM=\R^n$, for any $x\in M$,
and similarly identify the spaces of covectors,  $T_x^*M=(\R^n)^*$.

As earlier, we consider an initialized system
$\Sigma:$\ \  $x(t)=f(x(t-1),u(t)), \ \ x(0)=x_0$.
We define the Hamiltonian $H:M\times(\R^n)^*\times U\to \R$ of $\Sigma$,
\[
H(x,p,u)=pf(x,u),
\]
where we treat $p$ as a linear function acting on  $v=f(x,u)$, in coordinates $pv=\sum_ip_iv^i$.

\begin{remark}
If $M$ is an open subset of an affine space then a more natural Hamiltonian is
\[
H(x,p,u)=p(f(x,u)-x),
\]
since $v=f(x,u)-x$ is a vector and $pv$ is well defined by duality of covectors and vectors.
Further considerations hold with both definitions, with necessary modifications stated in remarks.
\end{remark}

In order to state the second order optimality conditions in terms of the Hamiltonian we consider a control sequence
$\bu=(u_1,\dots,u_N)$, where $u_t=u(t)$, $t=1,\dots,N$, and the corresponding trajectory $(x_0,\dots,x_N)$, $x_t=x(t)$. We can restrict the control sequence to its initial part
\[
\bu^i=(u_1,\dots,u_i),
\]
for any fixed $i=1,\dots,N$. All the definitions from the previous sections work if the sequence $\bu$ is replaced with
the restricted sequence $\bu^i$. In particular, using definitions from Section \ref{s-Local controllability}, we denote the vector fields and the Hessian matrix corresponding to the restricted sequence by
\[
L_\bu^i(x_0):=L_{\bu^i}(x_0), \qquad K_\bu^i(x_0):=K_{\bu^i}(x_0), \qquad  H_\bu^i:=H_{\bu^i}.
\]
Then
\[
L_\bu^1(x_0)\subset \cdots \subset L_\bu^N(x_0)
\]
and
\[
K_\bu^1(x_0)\subset \cdots \subset K_\bu^N(x_0)
\]
where, in the latter case, we use the natural embeddings
\[
\R^m\subset \R^m\times\R^m\subset \cdots\cdots \subset \R^m\times\stackrel{\mathrm{N\ times}}{\cdots}\times\R^m\
\]
given by $(u_1,\dots,u_i)\mapsto (u_1,\dots,u_i,0,\dots,0)$.

Consider again the optimality problem (P1) from Section \ref{s-Optimal control}, for system $\Sigma$ satisfying (A1)-(A4). Theorem \ref{t4} can be reformulated in terms of the Hamiltonian in the following way.

\begin{thm}\label{t6}
If $(u(1),\dots,u(N))\in (\mathrm{int}\,U)^N$ is an optimal control sequence for problem (P1) and  $x(0),\dots,x(N)$ is the corresponding trajectory, then there exists a sequence of covectors $p(0),\dots,p(N)$, with $p(N)=d\varphi(x_N)$, such that the following conditions are satisfied. The state equations hold
\[
\qquad\qquad x(t)=\frac{\partial H}{\partial p}(x(t-1),p(t),u(t)),  \quad\  t=1,\dots,N, \eqno{(\Sigma)\ }
\]
together with the adjoint equations
\[
\qquad\ p(t-1)=\frac{\partial H}{\partial x}(x(t-1),p(t),u(t)), \quad\    t=1,\dots,N, \eqno{(\Sigma^*)}
\]
the criticality condition
\[
\qquad\qquad \frac{\partial H}{\partial u}(x(t-1),p(t),u(t))=0, \qquad\ t=1,\dots,N, \eqno{(CC)}
\]
and the second order necessary condition
\[
\qquad\qquad\qquad  p(0) H^t_\bu\vert_{K^t_\bu(x_0)}\ \ge\ 0, \qquad\qquad t=1,\dots,N. \eqno{(SO)}
\]
\end{thm}

\begin{remark}\label{re-state-adjoint-eqns}
Note that equations $(\Sigma)$ together with $(\Sigma^*)$ are equivalent to
\[
\qquad\qquad x(t)=f(x(t-1),u(t)),  \qquad\qquad\ \ t=1,\dots,N,\  \eqno{(\Sigma)\ }
\]
\[
\qquad p(t-1)=p(t)\frac{\partial f}{\partial x}(x(t-1),u(t)), \qquad t=1,\dots,N. \eqno{(\Sigma^*)}
\]
\end{remark}

\begin{remark}
In the case of $M$ being a subset of an affine space and the Hamiltonian of the form $H(x,p,u)=p(f(x,u)-x)$ the theorem holds with the modified state equations
\[
 x(t)-x(t-1)=\frac{\partial H}{\partial p}(x(t-1),u(t)),  \quad\quad\ \ t=1,\dots,N, \eqno{(\Sigma)\ }
\]
and the adjoint equations
\[
\ p(t-1)-p(t)=\frac{\partial H}{\partial x}(x(t-1),p(t),u(t)), \quad   t=1,\dots,N, \eqno{(\Sigma^*)}
\]
\end{remark}

{\bf Proof of Theorem \ref{t6}.} Let $u(1)=u_1,\dots,u(N)=u_N$ be an optimal control sequence and  $x(1)=x_1,\dots,x(N)=x_N$ the corresponding trajectory starting from $x(0)=x_0$. By Theorem \ref{t4} the covector
\[
\lambda=d\varphi(x_N)\,df_{u_N}(x_{N-1})\cdots\cdots df_{u_1}(x_0)
\]
satisfies $\lambda\, Y^i_\bu(x_0)=0$ for $i=1,\dots,N$ (equivalently, $\lambda\in (L_\bu(x_0))^\perp)$) and
$\lambda\, H_\bu\vert_{K_\bu(x_0)}\ \ge \ 0$.

Denoting $p(i)=p_i$ we put
\[
p_N=\lambda\, (df_{u_1}(x_0))^{-1}\cdots (df_{u_N}(x_{N-1}))^{-1}
\]
and consider the sequence
\[
p_t=p_N\,df_{u_N}(x_{N-1})\cdots df_{u_{t+1}}(x_{t})=\lambda\, (df_{u_1}(x_0))^{-1}\cdots (df_{u_t}(x_{t-1}))^{-1},
\]
$t=1,\dots,N-1$ (for $t=0$ the latter equation does not hold and will not be used). Note that
\[
p_0=\lambda\quad  \mathrm{and} \quad p_N=d\varphi(x_N).
\]
The above sequence satisfies  $p_{t-1}=p_tdf_{u_t}(x_{t-1})=\partial H/\partial x(x_{t-1},p_t,u_t)$ which is the adjoint equation ($\Sigma^*$).  We claim that this solution satisfies the other assertions of the theorem, too.
The equation ($\Sigma$) is satisfied by the definitions of the sequences $u(t)$, $x(t)$ and $p(t)$ and Remark \ref{re-state-adjoint-eqns}.

We will show that (CC) follows from condition $\lambda\, Y^i_\bu(x_0)=0$. Namely, for $i=1$,
\begin{eqnarray*}
0 & = & \lambda\, Y^1_\bu(x_0)=\lambda\,X^+_{u_1}(x_0) = p_0\,(df_{u_1}(x_0))^{-1}\frac{\partial f_{u}}{\partial u}(x_0)\vert_{u=u_1}\\
  & = & p_1\frac{\partial f_u}{\partial u}(x_0)\vert_{u=u_1} = \frac{\partial H}{\partial u}(x_0,p_1,u_1).
\end{eqnarray*}
For $i=2,\dots,N$ we have
\begin{eqnarray*}
0 & = & \lambda\, Y^i_\bu(x_0)=\lambda\,(\Ad_{u_1}\cdots\Ad_{u_{i-1}}X^+_{u_i})(x_0)\\
  & = & \lambda\, (df_{u_1}(x_0))^{-1}\cdots (df_{u_{i-1}}(x_{i-2}))^{-1}X^+_{u_i}(x_{i-1})\\
  & = & p_{i-1}X^+_{u_i}(x_{i-1}) = p_{i-1}(df_{u_i}(x_{i-1}))^{-1}\frac{\partial f_{u}}{\partial u}(x_{i-1})\vert_{u=u_i}\\
  & = & p_i\frac{\partial f_{u}}{\partial u}(x_{i-1})\vert_{u=u_i} = \frac{\partial H}{\partial u_i}(x_{i-1},p_i,u_i),
\end{eqnarray*}
and we see that condition (CC) is satisfied.

In order to show that condition (SO) holds note first that it holds for $t=N$, by Theorem \ref{t4}. Namely, with our construction of the adjoint sequence $p(N),\dots,p(0)$ we have  $p_0=p(0)=\lambda$, where $\lambda$ is as in Theorem \ref{t4}.
We claim that for $1\le t<N$ condition (SO) follows from (SO) satisfied for $t=N$.
This is rather obvious. Namely, under the natural embedding $\R^{im}\subset \R^{Nm}$ defined by $(u_1,\dots,u_i)\mapsto (u_1,\dots,u_i,0\dots,0)$,
we have $K^i_\bu(x_0)\subset K^N_\bu(x_0)$ which follows from the definition of $K^i_\bu(x_0)$.
The matrix $H^i_\bu$ forms the upper-left block of $H^N_\bu$, of size $im$, corresponding to the subspace $\R^{im}\subset \R^{Nm}$. Thus positive definiteness of the whole matrix implies the same property of the sub-matrix.
The proof is complete.
\qed

Second order necessary conditions for geometric optimality can also be stated using Hamiltonian.

\begin{thm}\label{t7}
If $\bu=(u(1),\dots,u(N))\in (\mathrm{int}\, U)^N$ is a geometrically optimal control sequence of an initialized invertible system $\Sigma$ and  $x(1),\dots,x(N)$ is the corresponding trajectory starting from $x(0)=x_0$, then $\dim L_\bu(x_0)<n$ and, for any nonzero $\lambda\in (L_\bu(x_0))^\perp$, there exists a sequence of nonzero covectors $p(0),\dots,p(N)$ with $p(0)=\lambda$ such that the state equations ($\Sigma$), the adjoint equations ($\Sigma^*$), the criticality condition (CC) and the following index condition
\[
\qquad\quad  Ind^-(\lambda\, H^t_\bu)\vert_{K^t_\bu(x_0)} < \codim L^t_\bu(x_0), \qquad t=1,\dots,N, \eqno{(IC)}
\]
hold. In particular, (IC) implies that if $\mathrm{codim}\, L^t_\bu(x_0)=0$ for some $t\in \{1,\dots,N\}$ then the trajectory is not geometrically optimal. If $\mathrm{codim}\, L^t_\bu(x_0)=1$, then (IC) means that the quadratic form $(\lambda\, H^t_\bu)\vert_{K^t_\bu(x_0)}$ is non-negative definite.
\end{thm}

{\bf Proof of Theorem \ref{t7}.}  The proof is analogous to the proof of Theorem \ref{t6} with the difference that one should use Theorem \ref{t3}, instead of Theorem \ref{t4}, for determining the covector $\lambda$ with required properties. We leave the details to the reader.
\qed

\bibliographystyle{plain}
\bibliography{References}

\end{document}